\documentclass[10pt]{article}

\usepackage[latin1]{inputenc}
\usepackage{amsmath, amssymb, amsfonts, amsthm}
\usepackage{upgreek}
\usepackage{amsthm}
\usepackage{fullpage}
\usepackage{graphicx}
\usepackage{cancel}
\usepackage{subfigure}
\usepackage{mathrsfs}
\usepackage{outlines}
\usepackage[font={sf,it}, labelfont={sf,bf}, labelsep=space, belowskip=5pt]{caption}
\usepackage{hyperref}
\usepackage{titling}
\usepackage{xifthen}
\usepackage{color}

\usepackage{fancyhdr}
\usepackage[title]{appendix}
\usepackage{float}

\usepackage{bm}

\newcommand{\thedate}{August 8, 2019}
\date{\thedate}

\newcommand{\documenttitle}{Analytic Eigensystems for Isotropic Membrane Energies}

\DeclareMathOperator{\tr}{tr}

\newcommand{\defeq}{\vcentcolon=}
 
\allowdisplaybreaks

\providecommand{\norm}[1]{\lVert#1\rVert}

\providecommand{\pder}[2]{\frac{\partial #1}{\partial #2}}

\renewcommand{\vec}[1]{{\bf #1}}

\providecommand{\spderF}[1]{\frac{\partial^2 #1}{\partial F^2}}

\usepackage{prettyref}
\newrefformat{sec}{Section~\ref{#1}}
\newrefformat{tbl}{Table~\ref{#1}}
\newrefformat{fig}{Figure~\ref{#1}}
\newrefformat{chp}{Chapter~\ref{#1}}
\newrefformat{eqn}{\eqref{#1}}
\newrefformat{set}{\eqref{#1}}
\newrefformat{alg}{Algorithm~\ref{#1}}
\newrefformat{apx}{Appendix~\ref{#1}}
\newrefformat{prop}{Proposition~\ref{#1}}
\newcommand\pr[1]{\prettyref{#1}}

\def\normal{\hat {\bm n}}
\def\n{\normal}

\def\v{\vec{v}}

\def\R{\, \mathbb{R}}

\def\Ione{I_1^{3 \times 2}}
\def\Itwo{I_2^{3 \times 2}}
\def\Ithree{I_3^{3 \times 2}}

\newcommand*{\rom}[1]{\expandafter\@slowromancap\romannumeral #1@}
\newcommand{\RN}[1]{\textup{\uppercase\expandafter{\romannumeral#1}}}

\makeatletter
\usepackage{mathtools}
\newcases{mycases}{\quad}{%
  \hfil$\m@th\displaystyle{##}$}{$\m@th\displaystyle{##}$\hfil}{\lbrace}{.}
\makeatother

\title{\documenttitle}
\author{Julian Panetta}

\begin{document}
\maketitle

This document follows the approach of \cite{Smith2019} to derive the
Hessian eigenvalues and eigenmatrices for isotropic membrane energy densities $\psi(F)$, where
$F$ is a $3 \times 2$ deformation gradient. We assume that the energy is
expressed in terms of the following generalizations for $3 \times 2$ matrices
of the $2 \times 2$ tensor invariants\footnote{The $I_2$ invariant used here is
from \cite{Smith2019}; the
other standard definition of principal invariant $I_2 =
\frac{1}{2}\big(\tr(A)^2 - \norm{A}_F^2\big)$ actually coincides with $I_3$ in the 2D
case}:
\begin{align*}
    \Ione &\defeq \sigma_1 + \sigma_2 \\
    \Itwo &\defeq F : F = \sigma_1^2 + \sigma_2^2 \\
    \Ithree&\defeq \sigma_1 \sigma_2.
\end{align*}
In these definitions, $\sigma_1$ and $\sigma_2$ are the singular values of $F$
obtained from the singular value decomposition:
$$
F = U \underbrace{\begin{bmatrix} \sigma_1 & 0 \\ 0 & \sigma_2 \\ 0 & 0  \end{bmatrix}}_{\Sigma} V^T \quad \quad U \in O(3), V \in O(2).
$$
We note that the third column of $U$ is the deformed surface normal $\normal$.

\section{Differentiating the SVD}
We will need formulas for how $U$, $\Sigma$, and $V$ change as $F$ is perturbed with ``velocity''
$\dot F$, which we find by differentiating both sides of the SVD:
\begin{equation}
\label{eqn:svd_der_preliminary}
\dot F = \dot U \Sigma V^T + U \dot \Sigma V^T + U \Sigma \dot V^T
       \quad \Longrightarrow \quad U^T \dot F V = U^T \dot{U} \Sigma + \dot \Sigma + \Sigma \dot{V}^T V.
\end{equation}
Differentiating the relationships $U^T U = \text{Id}_{3 \times 3}$ and $V^T V =
\text{Id}_{2 \times 2}$ reveals that $U^T \dot{U}$ and $\dot{V}^T V$ are skew
symmetric and can be written as the infinitesimal rotations:
$$
U^T \dot{U} =
\begin{bmatrix}
    0 & -\omega_z & \omega_y \\
    \omega_z & 0 & -\omega_x \\
   -\omega_y & \omega_x & 0 \\
\end{bmatrix},
\quad \quad
\dot{V}^T V = 
\begin{bmatrix}
    0      & -\alpha \\
    \alpha & 0       \\
\end{bmatrix}.
$$
Plugging these into \pr{eqn:svd_der_preliminary}, we obtain a formula for
the infinitesimal rotations and singular value perturbations induced by $\dot{F}$:
\begin{equation}
\label{eqn:svd_der}
U^T \dot F V =
    \begin{bmatrix} \dot{\sigma_1} & -(\sigma_2 \omega_z + \sigma_1 \alpha) \\
                    \sigma_1 \omega_z + \sigma_2 \alpha & \dot{\sigma_2} \\
                   -\sigma_1 \omega_y & \sigma_2 \omega_x
    \end{bmatrix}.
\end{equation}
Geometrically, $\omega_z$ indicates a rotation of the surface element about the current normal $\n$, while
$\omega_x$ and $\omega_y$ are rotations around the principal stretch axes.
When $\omega_x = \omega_y = 0$, the deformed surface element simply rotates in-plane around $\n$
(and $\n$ does not change). However, nonzero $\omega_x$ and $\omega_y$ indicate
that $\dot{F}$ induces a rotation of $\n$.

\subsection{Example Perturbations}
\label{sec:example_perturb}
According to \pr{eqn:svd_der}, a perturbation of the form
$$
\dot{F} = U \begin{bmatrix} a & b \\ c & d \\ 0 & 0 \end{bmatrix} V^T
$$
leaves $\n$ unchanged as it stretches/rotates the surface element in-plane.
Specifically, we have $\dot{\sigma_1} = a$, $\dot{\sigma_2} = d$ and the following
system for $\omega_z$ and $\alpha$:
\begin{equation}
\label{eqn:rot2d_sum}
\begin{aligned}
    \sigma_2 \omega_z + \sigma_1 \alpha &= -b \\
    \sigma_1 \omega_z + \sigma_2 \alpha &=  c
\end{aligned}
\end{equation}
On the other hand, perturbation
$$
\dot{F} = U \begin{bmatrix} 0 & 0 \\ 0 & 0 \\ e & f \end{bmatrix} V^T
$$
rotates the surface element's normal by angular velocities $\omega_x = f / \sigma_2,
\omega_y = -e / \sigma_1$ without any in-plane stretch/rotation.

\section{Gradients of the Invariants}
We can now use the formulas for $\dot{\sigma_1}$ and $\dot{\sigma_2}$ to
differentiate the invariants:
\begin{align*}
\pder{\Ione}{F} : \dot{F} = \dot{\Sigma} : \begin{bmatrix} \ 1 \  & 0 \ \\ 0 & 1 \ \\ 0 & 0 \ \end{bmatrix} = \Big(U^T \dot{F} V\Big) : \begin{bmatrix} \ 1 \  & 0 \ \\ 0 & 1 \ \\ 0 & 0 \ \end{bmatrix}
    = \dot{F} : \left(U \begin{bmatrix} \ 1 \  & 0 \ \\ 0 & 1 \ \\ 0 & 0 \ \end{bmatrix} V^T\right) \quad \Longrightarrow \quad \pder{\Ione}{F} &= U \begin{bmatrix} \ 1 \  & 0 \ \\ 0 & 1 \ \\ 0 & 0 \ \end{bmatrix} V^T,
\\
\pder{\Ithree}{F} : \dot{F} = \dot{\Sigma} : \begin{bmatrix} \sigma_2 & 0 \\ 0 & \sigma_1 \\ 0 & 0 \end{bmatrix} = \Big(U^T \dot{F} V\Big) : \begin{bmatrix} \sigma_2 & 0 \\ 0 & \sigma_1 \\ 0 & 0 \end{bmatrix}
    = \dot{F} : \left(U \begin{bmatrix} \sigma_2 & 0 \\ 0 & \sigma_1 \\ 0 & 0 \end{bmatrix} V^T\right) \quad \Longrightarrow \quad
   \pder{\Ithree}{F} &= U \begin{bmatrix} \sigma_2 & 0 \\ 0 & \sigma_1 \\ 0 & 0 \end{bmatrix} V^T,
\\
    \pder{\Itwo}{F} : \dot{F} = 2 F : \dot{F} \quad \Longrightarrow \quad \pder{\Itwo}{F} &= 2 F.
\end{align*}

\section{Hessians of the Invariants}
We evaluate the Hessian applied to an
arbitrary perturbation $\dot{F}$. First, the easy invariant:
$$
\spderF{\Itwo} : \dot{F} = 2 \dot{F},
$$
which means $\spderF{\Itwo}$ is a multiple of the fourth order identity
tensor. Any orthogonal basis can be chosen as a set of eigenmatrices, and their
corresponding eigenvalues are all $2$.

Next, we consider $\Ione$:
$$
U^T \left( \spderF{\Ione} : \dot{F} \right) V
    = U^T \dot{U} \begin{bmatrix} 1 & 0 \\ 0 & 1 \\ 0 & 0 \end{bmatrix} +
        \begin{bmatrix} 1 & 0 \\ 0 & 1 \\ 0 & 0 \end{bmatrix} \dot{V}^T V
            = \begin{bmatrix} 0 & -(\omega_z + \alpha) \\ \omega_z + \alpha & 0 \\ -\omega_y & \omega_x \end{bmatrix}.
$$
We plug in $\dot{F} = U \begin{bmatrix} a & b \\ c & d \\ e & f\end{bmatrix} V^T$
and note that summing the equations in \pr{eqn:rot2d_sum} yields $\omega_z + \alpha = \frac{c - b}{\sigma_1 + \sigma_2}$.
Thus:
$$
\spderF{\Ione} : \dot{F}
    = U \begin{bmatrix} 0 & \frac{b - c}{\sigma_1 + \sigma_2} \\
                        \frac{c - b}{\sigma_1 + \sigma_2} & 0 \\
                        \frac{e}{\sigma_1} & \frac{f}{\sigma_2} \end{bmatrix} V^T.
$$
From this expression, we see there is a three dimensional null space
with $e = f = 0$ and $b = c$. We can pick the following
orthonormal basis for this subspace:
$$
\frac{1}{\sqrt{2}} U \begin{bmatrix} 1 & 0 \\ 0 &  1 \\ 0 & 0 \end{bmatrix} V^T,\quad
\frac{1}{\sqrt{2}} U \begin{bmatrix} 1 & 0 \\ 0 & -1 \\ 0 & 0 \end{bmatrix} V^T,\quad
\frac{1}{\sqrt{2}} U \begin{bmatrix} 0 & 1 \\ 1 &  0 \\ 0 & 0 \end{bmatrix} V^T \quad \quad (\lambda = 0).
$$
We further deduce the three eigenmatrices with nonzero eigenvalues:
$$
\underbrace{\frac{1}{\sqrt{2}}U \begin{bmatrix} 0 & -1 \\ 1 &  0 \\ 0 & 0 \end{bmatrix} V^T}_{\lambda = \frac{2}{\sigma_1 + \sigma_2}}, \quad
\underbrace{                  U \begin{bmatrix} 0 &  0 \\ 0 &  0 \\ 1 & 0 \end{bmatrix} V^T}_{\lambda = \frac{1}{\sigma_1}}, \quad
\underbrace{                  U \begin{bmatrix} 0 &  0 \\ 0 &  0 \\ 0 & 1 \end{bmatrix} V^T}_{\lambda = \frac{1}{\sigma_2}}.
$$

Finally, we consider $\Ithree$:
$$
U^T \left( \spderF{\Ithree} : \dot{F} \right) V
    = U^T \dot{U} \begin{bmatrix} \sigma_2 & 0 \\ 0 & \sigma_1 \\ 0 & 0 \end{bmatrix}
        + 
        \begin{bmatrix} \dot \sigma_2 & 0 \\ 0 & \dot \sigma_1 \\ 0 & 0 \end{bmatrix}
        +
        \begin{bmatrix} \sigma_2 & 0 \\ 0 & \sigma_1 \\ 0 & 0 \end{bmatrix} \dot{V}^T V
    = \begin{bmatrix} \dot \sigma_2 & -(\sigma_1 \omega_z + \sigma_2 \alpha) \\
                      \sigma_2 \omega_z + \sigma_1 \alpha & \dot \sigma_1 \\
                     -\sigma_2 \omega_y & \sigma_1 \omega_x \end{bmatrix}.
$$
Again plugging in $\dot{F} = U \begin{bmatrix} a & b \\ c & d \\ e & f\end{bmatrix} V^T$ and
using the formulas from \pr{sec:example_perturb}, we find:
$$
\spderF{\Ithree} : \dot{F} =
U \begin{bmatrix}  d & -c \\
                  -b & a \\
\frac{\sigma_2}{\sigma_1} e & \frac{\sigma_1}{\sigma_2} f \end{bmatrix} V^T
$$
We deduce the following eigenmatrices and eigenvalues:
\begin{gather*}
    \underbrace{\frac{1}{\sqrt{2}} U \begin{bmatrix} 1 &  0 \\ 0 &  1 \\ 0 & 0 \end{bmatrix} V^T,\quad
                \frac{1}{\sqrt{2}} U \begin{bmatrix} 0 & -1 \\ 1 &  0 \\ 0 & 0 \end{bmatrix} V^T}_{\lambda = 1},\quad 
    \underbrace{\frac{1}{\sqrt{2}} U \begin{bmatrix} 1 &  0 \\ 0 & -1 \\ 0 & 0 \end{bmatrix} V^T,\quad
                \frac{1}{\sqrt{2}} U \begin{bmatrix} 0 &  1 \\ 1 &  0 \\ 0 & 0 \end{bmatrix} V^T}_{\lambda = -1},\quad \\
    \underbrace{                   U \begin{bmatrix} 0 &  0 \\ 0 &  0 \\ 1 & 0 \end{bmatrix} V^T}_{\lambda = \frac{\sigma_2}{\sigma_1}},\quad
    \underbrace{                   U \begin{bmatrix} 0 &  0 \\ 0 &  0 \\ 0 & 1 \end{bmatrix} V^T}_{\lambda = \frac{\sigma_1}{\sigma_2}}.
\end{gather*}

We note that for all invariants, four of the six Hessian eigenmatrices are
simply padded versions of the 2D eigenmatrices from \cite{Smith2019}, while the
last two are new and concern the rotation of the surface element's normal.

\section{Example: Incompressible neo-Hookean Sheet}
We consider the membrane energy of a thin sheet of incompressible
neo-Hookean material \cite{bonet1997nonlinear}:
$$
\psi_\text{IncNeo}(F_\text{3D}) = \frac{\mu}{2} \left(\tr(F_\text{3D}^T F_\text{3D}) - 3\right)
                                = \frac{\mu}{2} \left(I_2^\text{3D} - 3\right)
$$
When the sheet experiences an in-plane deformation gradient $F \in \R^{3 \times 2}$,
it stretches or compresses in the normal direction to maintain $J = 1$.
We can solve for the normal stretch as $\frac{1}{\Ithree}$ and express
$\psi_\text{IncNeo}$ directly in terms of $F$'s invariants:
$$
\psi_\text{sheet}(F) = \frac{\mu}{2} \left( \Itwo + \left(\frac{1}{\Ithree}\right)^2 - 3 \right).
$$
The Hessian of this energy density is:
\begin{align*}
\spderF{\psi_\text{sheet}} &=
    \frac{\mu}{2} \left[
        \spderF{\Itwo} + 6 \left(\frac{1}{\Ithree}\right)^4 \pder{\Itwo}{F} \otimes \pder{\Itwo}{F}
                       - 2 \left(\frac{1}{\Ithree}\right)^3 \spderF{\Ithree}
    \right]
\\
    &=
    \mu \left[
        \text{Id}_4 + 3 \left(\frac{1}{\Ithree}\right)^4 \left(U \begin{bmatrix} \sigma_2 & 0 \\ 0 & \sigma_1 \\ 0 & 0 \end{bmatrix} V^T\right) \otimes \left(U \begin{bmatrix} \sigma_2 & 0 \\ 0 & \sigma_1 \\ 0 & 0 \end{bmatrix} V^T\right)
                    - \left(\frac{1}{\Ithree}\right)^3 \spderF{\Ithree}
    \right].
\end{align*}
We note that $\pder{\Itwo}{F}$ is orthogonal to all but two of the eigenmatrices of $\spderF{\Ithree}$
(and eigenmatrices for the fourth order identity tensor $\text{Id}_4$ can be chosen arbitrarily), so we
immediately get the following four eigenpairs:
$$
    \underbrace{\frac{1}{\sqrt{2}} U \begin{bmatrix} 0 & -1 \\ 1 &  0 \\ 0 & 0 \end{bmatrix} V^T}_{\lambda = \mu  - \mu \left(\frac{1}{\Ithree}\right)^3},\quad 
    \underbrace{\frac{1}{\sqrt{2}} U \begin{bmatrix} 0 &  1 \\ 1 &  0 \\ 0 & 0 \end{bmatrix} V^T}_{\lambda = \mu  + \mu \left(\frac{1}{\Ithree}\right)^3},\quad
    \underbrace{                   U \begin{bmatrix} 0 &  0 \\ 0 &  0 \\ 1 & 0 \end{bmatrix} V^T}_{\lambda = \mu  - \mu \left(\frac{1}{\Ithree}\right)^3\frac{\sigma_2}{\sigma_1}},\quad
    \underbrace{                   U \begin{bmatrix} 0 &  0 \\ 0 &  0 \\ 0 & 1 \end{bmatrix} V^T}_{\lambda = \mu  - \mu \left(\frac{1}{\Ithree}\right)^3\frac{\sigma_1}{\sigma_2}}.
$$
Because $\pder{\Itwo}{F}$ is generally not orthogonal to either of the remaining two
eigenmatrices of $\spderF{\Ithree}$ (whose eigenvalues are distinct) we must diagonalize
the projection of $\spderF{\psi_\text{sheet}}$ onto their span to obtain the final two
eigenpairs. We obtain simpler expressions using the
basis $D_1 \defeq U \begin{bmatrix} 1 & 0 \\ 0 &  0 \\ 0 & 0 \end{bmatrix} V^T$
and $D_2 \defeq U \begin{bmatrix} 0 & 0 \\ 0 &  1 \\ 0 & 0 \end{bmatrix} V^T$
for this subspace, which results in the reduced Hessian:
\begin{align*}
&\begin{bmatrix}
    D_1 : \spderF{\psi_\text{sheet}} : D_1 &
    D_1 : \spderF{\psi_\text{sheet}} : D_2 \\
    D_2 : \spderF{\psi_\text{sheet}} : D_1 &
    D_2 : \spderF{\psi_\text{sheet}} : D_2
\end{bmatrix}
% \\ &\quad
= \mu
\begin{bmatrix}
    1 & 0 \\
    0 & 1
\end{bmatrix}
+
\frac{\mu}{\left(\Ithree\right)^4}
\begin{bmatrix}
    3 \sigma_2^2 & 2 \Ithree \\
       2 \Ithree & 3 \sigma_1^2
\end{bmatrix}.
\end{align*}
The eigendecomposition of this $2 \times 2$ matrix can be expressed
by introducing quantities ${\beta \defeq 3 (\sigma_2^2 - \sigma_1^2)}$ and
${\gamma \defeq \sqrt{16 \left(\Ithree\right)^2 + \beta^2}}$:
$$
\v_1 = \begin{bmatrix}\beta - \gamma \\ 4 \Ithree \end{bmatrix},\quad \lambda_1 = \mu  + \mu \frac{3 \Itwo + \gamma}{2 \left(\Ithree\right)^4}, \quad \quad \quad \quad
\v_2 = \begin{bmatrix}\beta + \gamma \\ 4 \Ithree \end{bmatrix},\quad \lambda_2 = \mu  + \mu \frac{3 \Itwo + \gamma}{2 \left(\Ithree\right)^4},
$$
making the final two eigenpairs of $\spderF{\psi_\text{sheet}}$:
$$
    \underbrace{U \begin{bmatrix} \beta - \gamma & 0 \\ 0 &  4 \Ithree \\ 0 & 0 \end{bmatrix} V^T}_{\lambda = \mu  + \mu \frac{3 \Itwo + \gamma}{2 \left(\Ithree\right)^4}},\quad 
    \underbrace{U \begin{bmatrix} \beta + \gamma & 0 \\ 0 &  4 \Ithree \\ 0 & 0 \end{bmatrix} V^T}_{\lambda = \mu  + \mu \frac{3 \Itwo - \gamma}{2 \left(\Ithree\right)^4}}.
$$
Note that these eigenmatrices do not have unit norm and should be normalized.

\bibliographystyle{plain}
\bibliography{AnalyticEigenvalueIsotropicMembranes}

\end{document}